\documentclass[12pt]{amsart}
\usepackage{amsmath,amscd,amssymb,amsfonts}
\def\vers{Aug.~8, 2010, v.3b}
\def\1{\hskip1pt}
\def\scirc{\,\raise.2ex\h{${\scriptstyle\circ}$}\,}
\def\ssb{\raise.2ex\h{${\scriptscriptstyle\bullet}$}}
\def\don{\rlap{$\downarrow$}\raise2pt\hbox{$\downarrow$}}
\def\msum{\h{$\sum$}}
\def\mpro{\h{$\prod$}}
\def\mcap{\h{$\bigcap$}}

\def\mopl{\h{$\bigoplus$}}
\def\bs{\bigskip}
\def\ms{\medskip}
\def\ssk{\smallskip}
\def\nin{\noindent}
\def\A{{\mathcal A}}
\def\a{\alpha}
\def\C{{\mathbf C}}
\def\D{{\mathbf D}}
\def\DD{{\mathcal D}}
\def\d{{\mathbf d}}
\def\da{\downarrow}
\def\E{{{'\!E}}}
\def\F{{\mathcal F}}
\def\G{{\mathcal G}}
\def\g{\gamma}
\def\H{{\mathcal H}}
\def\HH{{}^p{\mathcal H}}
\def\tH{\widetilde{H}}
\def\h{\hbox}
\def\I{{\mathcal I}}
\def\i{\tilde{i}}

\def\OO{{\mathcal O}}
\def\P{{\mathbf P}}
\def\ph{\varphi}
\def\Q{{\mathbf Q}}
\def\q{\quad}

\def\R{{\mathbf R}}
\def\S{\Sigma}
\def\Y{{\mathcal Y}}
\def\Z{{\mathbf Z}}
\def\Gr{{\rm Gr}}
\def\IC{{\rm IC}}
\def\Coim{{\rm Coim}}
\def\Ext{{\rm Ext}}
\def\Hom{{\rm Hom}}
\def\Im{{\rm Im}}
\def\Ker{{\rm Ker}}
\def\MHM{{\rm MHM}}
\def\Perv{{\rm Perv}}
\def\Sing{{\rm Sing}\,}
\def\supp{{\rm supp}\,}
\def\prim{{\rm prim}}
\def\inv{{\rm inv}}
\def\van{{\rm van}}
\def\codim{{\rm codim}}

\def\into{\hookrightarrow}
\def\onto{\mathop{\rlap{$\to$}\hskip2pt\hbox{$\to$}}}
\def\simto{\buildrel\sim\over\to}
\def\({{\rm (}}
\def\){{\rm )}}

\begin{document}
\title[Vanishing cycle sheaves and quasi-semistable degenerations]
{Vanishing cycle sheaves of one-parameter smoothings
and quasi-semistable degenerations}
\author{Alexandru Dimca}
\address{Laboratoire J.A.\ Dieudonn\'e, UMR du CNRS 6621,
Universit\'e de Nice-Sophia Antipolis, Parc Valrose,
06108 Nice Cedex 02, France}
\email{Alexandru.DIMCA@unice.fr}
\author{Morihiko Saito}
\address{RIMS Kyoto University, Kyoto 606-8502 Japan}
\email{msaito@kurims.kyoto-u.ac.jp}
\thanks{This work is partially supported by Kakenhi 19540023 and by  ANR-08-BLAN-0317-02 (SEDIGA).
 Additionally, A.~Dimca is grateful to 
ASSMS, Government College University, Lahore, Pakistan, 
where part of the work on this paper was done.}
\date{\vers}
\begin{abstract}
We study the vanishing cycles of a one-parameter smoothing of a
complex analytic space and show that the weight filtration on its
perverse cohomology sheaf of the highest degree is quite close to
the monodromy filtration so that its graded pieces have a modified
Lefschetz decomposition. We describe its primitive part using the
weight filtration on the perverse cohomology sheaves of the
constant sheaves.
As a corollary we show in the local complete intersection case
that 1 is not an eigenvalue of the monodromy on the reduced
Milnor cohomology at any points if and only if the total space
and the singular fiber are both rational homology manifolds.
Also we introduce quasi-semistable degenerations and
calculate the limit mixed Hodge structure by constructing the
weight spectral sequence.
As a corollary we show non-triviality of the space of vanishing
cycles of the Lefschetz pencil associated with a tensor product
of any two very ample line bundles except for the case of
even-dimensional projective space where two has to be replaced
by three.
\end{abstract}

\maketitle

\centerline{\bf Introduction}

\bs\nin
Let $X$ be a complex analytic space of dimension $n$, and $X_i$ be
the irreducible components of dimension $n$.
Let $\IC_X\Q:=\mopl_i\IC_{X_i}\Q$ be the intersection complex,
which belongs to the category of perverse sheaves $\Perv(X,\Q)$
(see [3])
and is naturally identified with a mixed Hodge module ([23], 3.21).
Shrinking $X$ if necessary, we assume that $X$ is an intersection
of hypersurfaces in a complex manifold so that $\Q_X$ exists in
the derived category of mixed Hodge modules, see (1.3) below.
By the same argument as in [23], 4.5.9,
there is a canonical morphism
$$\Q_X[n]\to\IC_X\Q,$$
inducing an isomorphism
$$\Gr^W_n\HH^n\Q_X\simto\IC_X\Q,\leqno(0.1)$$
where $\HH^j$ is the perverse cohomology [3] and $W$ is the weight
filtration of $\HH^j\Q_X$ defined in the abelian category of mixed
Hodge modules $\MHM(X)$, see [23].
(In this paper mixed Hodge modules are denoted by their underlying
perverse sheaves when there is no fear of confusions.)
It is well-known (see [3], 4.2.4 and [23], 2.26) that
$$\HH^j\Q_X=0\,\,(j>n),\q\Gr^W_k\HH^j\Q_X=0\,\,(k>j).
\leqno(0.2)$$
It is also well-known to the specialists that $\Q_X[n]=\IC_X\Q$
if and only if $X$ is a $\Q$-homology manifold, see (1.5) below.
Note that the $\HH^j\Q_X$ are globally well-defined in $\MHM(X)$
(see [23], 2.19) and we have by (0.1--2)
$$\Q_X[n]=\IC_X\Q\iff\big[\1\HH^j\Q_X=0\,\,(j<n),\,\,
\Gr^W_j\HH^n\Q_X=0\,\,(j<n)\big].\leqno(0.3)$$

Let $f:X\to\C$ be a holomorphic function on $X$ such that the
fibers $f^{-1}(t)$ are smooth and purely $(n-1)$-dimensional
for $t\ne 0$.
Setting $Y=f^{-1}(0)$, this condition is essentially equivalent
to the smoothness and the pure dimensionality of
$X\setminus Y$ replacing $X$ with an open neighborhood of $Y$
if necessary.
Assume $f$ is nonconstant on any $n$-dimensional irreducible
component of $X$.
Let $i:Y\into X$ denote the inclusion.

\ms\nin
{\bf Theorem~1.} {\it Let $f$ be as above. Then we have an
exact sequence of mixed Hodge modules on $Y$
$$0\to\HH^{n-1}\Q_X\to\HH^{n-1}\Q_Y\to\HH^{-1}i^*(\IC_X\Q)\to
W_{n-1}\HH^n\Q_X\to 0,$$
and isomorphisms}
$$\HH^j\Q_X=\HH^j\Q_Y\q(j<n-1).$$

\ms
This is shown by applying the functor $i^*$ to the distinguished
triangle associated to the mapping cone of $\Q_X\to(\IC_X\Q)[-n]$
and using (0.1--2).
If $\Q_X[n]$ is a perverse sheaf, then $\Q_Y[n-1]$ is also a
perverse sheaf by the last isomorphisms of Theorem~1, and we get

\ms\nin
{\bf Corollary~1.} {\it If $\Q_X[n]$ is a perverse sheaf, then
we have $\IC_X\Q=\Q_X[n]$ if and only if
$i^*(\IC_X\Q)[-1]=\Q_Y[n-1]$.}

\ms
Let $\psi_f,\ph_f$ denote the nearby and vanishing cycle functors
[6].
We denote by $\psi_{f,1},\ph_{f,1}$ their unipotent monodromy part.
In this paper these functors are {\it shifted by} $-1$ so that they
preserve perverse sheaves and commute with $\HH^j$.
Set $N=\log T_u$ with $T_u$ the unipotent part of the monodromy $T$.
Note that $\psi_f\Q_X[n]=\psi_f\IC_X\Q$ is a perverse sheaf, and
$\psi_f\1\HH^j\Q_X=0$ for $j\ne n$,
since $X\setminus Y$ is smooth and $\psi_f\F$ depends only on
$\F|_{X\setminus Y}$ in general.
The weight filtration $W$ on the nearby cycles
$\psi_f\Q_X[n]=\psi_f\IC_X\Q$ is the monodromy filtration [7]
shifted by $n-1$.
However, $W$ on the vanishing cycles $\ph_{f,1}\1\HH^j\Q_X$ is
rather complicated in case $X$ is singular
(even in the isolated singularity case, see e.g.\ [27]).
Since $\HH^j\Q_X$ for $j\ne n$ is supported in $Y$,
$\ph_{f,1}\1\HH^j\Q_X$ is identified with $\HH^j\Q_X$ and the action
of $N$ on it vanishes.
So we are mainly interested in $\ph_{f,1}\1\HH^0(\Q_X[n])$.
Set
$$V_{\ssb}=\mopl_kV_k\q\h{with}\q V_k:=\Gr^W_k\ph_{f,1}\1\HH^0
(\Q_X[n])\in\MHM(Y).$$

\ms\nin
{\bf Theorem~2.} {\it Let $f$ be as in Theorem~$1$.
There is a noncanonical decomposition
$$\aligned&V_{\ssb}\cong V_{\ssb}'\oplus V_{\ssb}''\q\h{in}\,\,\,
\MHM(Y)\,\,\,\h{compatible with $N$ and satisfying}\\
&N^i:V_{n+i}'\simto V_{n-i}'(-i),\q
N^i:V_{n-1+i}''\simto V_{n-1-i}''(-i)\q\h{for}\,\,\,i>0,\endaligned$$
where $(-i)$ is the Tate twist {\rm [5]}.
Moreover, setting $K'_k=\Ker\,N\cap V'_k$ and
$K''_k=\Ker\,N\cap V''_k$, we have}
$$\aligned K'_k&=(\Gr_{k-2}^W(\HH^{n-1}\Q_Y/\HH^{n-1}\Q_X))(-1)\q
(k\le n),\\ K''_k&=\Gr_k^W\HH^n\Q_X\q(k\le n-1).\endaligned$$

\ms\nin
{\bf Corollary~2.} {\it With the above notation,
we have the Lefschetz decompositions
$$\aligned
&V'_{\ssb}=\mopl_{k\ge 0}\mopl_{i=0}^k\,N^iP'_{n+k}(i),\q
V''_{\ssb}=\mopl_{k\ge 0}\mopl_{i=0}^k\,N^iP''_{n-1+k}(i),\\
&\h{with}\q K'_{n-k}=N^kP'_{n+k}(k),\q
K''_{n-1-k}=N^kP''_{n-1+k}(k),\endaligned$$
where $P'_{n+k}:=\Ker\,N^{k+1}\subset V'_{n+k}$ denotes the
$N$-primitive part, and similarly for $P''_{n-1+k}$.}

\ms
If $\Q_X[n]$ is a perverse sheaf (e.g.\ if $X$ is a local complete
intersection), then $\HH^{n-1}\Q_X$ vanishes in the above formula
and $\Q_Y[n-1]$ is also a perverse sheaf by Theorem~1.
In the isolated singularity case, Theorem~2 was essentially stated
in [21], see also [27].

From Theorem~2 together with (0.3) for $X,Y$, we deduce

\ms\nin
{\bf Theorem~3.} {\it Let $f$ be as in Theorem~$1$.
Assume $\Q_X[n]$ is a perverse sheaf
\(e.g.\ $X$ is a local complete intersection\) so that
$\Q_Y[n-1]$ is also a perverse sheaf.
Let $F_x$ denote the Milnor fiber of $f$ around $x\in Y$.
Then

\ssk\nin
{\rm (a)}
The following three conditions are equivalent to each other.

\nin
\q{\rm (i)} $V'_{\ssb}=0$,

\nin
\q{\rm (ii)} $\Q_Y[n-1]=\IC_Y\Q$, i.e.\ $Y$ is a $\Q$-homology
manifold,

\nin
\q{\rm (iii)} $W$ on $\ph_{f,1}(\Q_X[n])$ is the monodromy
filtration shifted by $n-1$.

\ssk
\nin
{\rm (b)}
The following three conditions are equivalent to each other.

\nin
\q{\rm (i)} $V''_{\ssb}=0$,

\nin
\q{\rm (ii)} $\Q_X[n]=\IC_X\Q$, i.e.\ $X$ is a $\Q$-homology
manifold,

\nin
\q{\rm (iii)} $W$ on $\ph_{f,1}(\Q_X[n])$ is the monodromy
filtration shifted by $n$.

\ssk\nin
{\rm (c)}
The following three conditions are equivalent to each other.

\nin
\q{\rm (i)} $\ph_{f,1}\Q_X=0$,

\nin
\q{\rm (ii)} $X$ and $Y$ are $\Q$-homology manifolds,

\nin
\q{\rm (iii)} $1$ is not an eigenvalue of the monodromy on
$\tH^j(F_x,\Q)$ for any $j,x$.}

\ms
The assertion (c) says that $\ph_{f,1}\Q_X$ just contains the
information of the difference between $\Q_X[n]$ and $\IC_X\Q$
together with the difference between $\Q_Y[n-1]$ and $\IC_Y\Q$
in this case, see (0.3).
(For assertions using the topological methods, see [18].)
We have $V''_{\ssb}\ne 0$, for example, if $f:X\to\C$ is the
base change of $g:Z\to\C$ by an $m$-fold ramified covering of
$\C$ such that $Z$ is smooth and an $m$-th primitive root of
unity is an eigenvalue of the monodromy of $g$.
We have $V'_{\ssb}\ne 0$, for example, in the case where $X$ is
smooth, $Y$ has an isolated singularity, and $1$ is an eigenvalue
of the Milnor monodromy.

Let $z\in Z:=\supp V_{\ssb}\subset Y$ with the inclusion
$i_z:\{z\}\into Z$.
Let $Z_i$ be the local irreducible components of $(Z,z)$.
Using Theorem~2, we get an assertion on the Milnor cohomology
in a special case as follows.

\ms\nin
{\bf Corollary~3.} {\it Let $f$ be as in Theorem~$1$.
Assume $Z$ is a curve and all the eigenvalues of the monodromy
around $z$ of the local system $(V_{\ssb}|_{Z_i\setminus\{z\}})[-1]$
are different from $1$ for any $i$.
Then $H^ji_z^*V_{\ssb}=0$ for $j\ne 0$,
and the assertion of Theorem~$2$ holds by replacing
$V_{\ssb}$ with the graded pieces of the unipotent monodromy part
of the Milnor cohomology $\Gr^W_{\ssb}H^{n-1}(F_z,\Q)_1$, and
applying the functor $H^0i_z^*$ to the mixed Hodge modules
appearing in Theorem~$2$, where
$$H^0i_z^*K''_k=\Gr^W_kH^{-1}i_z^*\IC_X\Q.$$
In case $\HH^{n-1}\Q_X=0$, we have moreover
$$H^0i_z^*K'_k(1)=\Gr^W_{k-2}H^{-1}i_z^*\,C(\Q_Y[n-1]\to\IC_Y\Q).$$
Here the mapping cone $C(\Q_Y[n-1]\to\IC_Y\Q)$ may be replaced by
$\IC_Y\Q$ if $n\ge 3$.}

\ms\nin
{\bf Corollary~4.} {\it
With the assumptions of Corollary~$3$, the maximal size of the
Jordan block of the monodromy on $H^{n-1}(F_z,\Q)_1$ is the
largest number $k$ such that $H^0i_z^*K'_{n+1-k}\ne 0$ or
$H^0i_z^*K''_{n-k}\ne 0$ if $k>0$.}

\ms
As for the nearby cycles $\psi_{f,1}(\Q_X[n])$, we have the
following

\ms\nin
{\bf Theorem~4.} {\it Let $f$ be as in Theorem~$1$.
There are noncanonical decompositions for $k<n$}
$$\Gr^W_k\psi_{f,1}(\Q_X[n])\cap\Ker\,N\cong
\Gr_k^W(\HH^{n-1}\Q_Y/\HH^{n-1}\Q_X)\oplus\Gr_k^W\HH^n\Q_X.$$

\ms
Note that the weight filtration $W$ on $\psi_{f,1}(\Q_X[n])$ is
the monodromy filtration shifted by $n-1$, and the left-hand side
essentially gives the $N$-primitive part.

\ms
As an application of Theorem~4 we construct the weight spectral
sequence associated with a {\it quasi-semistable} degeneration,
see (3.1--2).
In a typical case, we have the following.

\ms\nin
{\bf Theorem~5.} {\it
Let $L_k\,(1\le k\le r)$ be line bundles on a smooth proper complex
algebraic variety $Y$ of dimension $n$, and $L_0$ be the tensor
product of $L_k\,(1\le k\le r)$ where $n,r\ge 2$.
Let $Y_k$ be smooth divisors defined by $g_k\in\Gamma(Y,L_k)\,
(0\le k\le r)$.
Assume $\bigcup_{k=0}^rY_k$ is a divisor with normal crossings on
$Y$.
Let $X=\{g_1\cdots g_r=tg_0\}\subset Y\times\C$ with $f:X\to\C$
induced by the second projection where $t$ is the coordinate of
$\C$.
Then we have a weight spectral sequence degenerating at $E_2$
$$E_1^{-k,j+k}=H^j(X_0,\Gr^W_{n-1+k}\psi_{f,1}\IC_X\Q)\Rightarrow
H^{j+n-1}(X_{\infty},\Q),$$
where $H^{\ssb}(X_{\infty},\Q)$ denotes the limit mixed Hodge
structure {\rm [25], [26]}.
Moreover, $E_1^{-i,j+i}$ is the direct sum of
$$\Bigl(\bigoplus_{|I|=i+2l+1}H^{j+n-|I|}(Y_I)(-i-l)\Bigr)\mopl
\Bigl(\bigoplus_{|I|=i+2l+2}H^{j+n-|I|-1}(Y'_I)(-i-l-1)\Bigr),$$
over $l\ge\max(-i,0)$, where
$Y'_I=Y_0\cap\bigcap_{i\in I}Y_i$ and
$Y_I=\bigcap_{i\in I}Y_i$ for $I\subset\{1,\dots,r\}$.}

\ms
The first part of the $E_1$-term is the same as the weight
spectral sequence in case of a semistable degeneration [26]
where $X$ is smooth and $X_0$ is a reduced divisor with simple
normal crossings.
The second part of the $E_1$-term is closely related to the
singularities of $X$.
Note that $X$ is {\it singular} and $f$ is {\it not} semistable
if $n\ge 3$.
Nevertheless we get the weight spectral sequence as above
{\it without} using a blow-up of $X$ to get a (non-reduced)
semistable model as in [11], [13].

In case $Y=\P^n$ we can explicitly describe the
$E_2$-term of the weight spectral sequence, i.e. the graded pieces
of the weight filtration of the limit mixed Hodge structure,
by using the $N$-primitive decomposition as below.

\ms\nin
{\bf Corollary~5.} {\it With the notation and the assumptions of
Theorem~$5$, assume $Y=\P^n$.
Let $H^{n-1}_{\prim}(X_{\infty},\Q)$, $H^{n-|I|}_{\prim}(Y_I,\Q)$,
$H^{n-1-|I|}_{\prim}(Y'_I,\Q)$ denote the middle primitive
cohomology, and $P_N\Gr^W_{n-1+k}H^{n-1}_{\prim}(X_{\infty},\Q)$
denote the $N$-primitive part defined by $\Ker\,N^{k+1}\,\,(k\ge 0)$.
Set $m=\binom{r-1}{n}$.

Then $P_N\Gr^W_{n-1+k}H^{n-1}_{\prim}(X_{\infty},\Q)$ for
$k\in[0,n-1]$ is given by}
$$\aligned&\bigl(\mopl_{|I|=k+1}\,\tH^{n-|I|}_{\prim}(Y_I)(-k)
\bigr)\h{$\oplus$}\bigl(\mopl_{|I|=k+2}\,\tH^{n-1-|I|}_{\prim}(Y'_I)
(-k-1))\,\,\,\!\h{if}\,\,k\ne n-1,\\
&\bigl(\mopl_{|I|=k+1}\,\tH^{n-|I|}_{\prim}(Y_I)(-k)\bigr)\h{$\oplus$}
\bigl(\buildrel{m}\over\mopl\Q(1-n)\bigr)\hskip3.2cm\!\h{if}\,\,
k=n-1.\endaligned$$

\ms
Here $\tH^{n-|I|}_{\prim}(Y_I)=0$ if $|I|>\min(n,r)$,
and $\tH^{n-1-|I|}_{\prim}(Y'_I)=0$ if $|I|>\min(n-1,r)$.
Let $d_k$ be the integers such that $L_k=\OO_{\P^n}(d_k)\,
(0\le k\le r)$.
Note that $d_0=\sum_{k=1}^rd_k$.
If $d_k=1 $ for any $r\ge 1$ as in [11], then $\tH^{n-|I|}_{\prim}
(Y_I)=0$ for any $I$, and $\Gr^W_{n-1+k}H^{n-1}(X_{\infty},\Q)$
has level $<n-1-|k|$ for $|k|\ne n-1$ (where the level means the
difference between the maximum and the minimum of the $p$ with
$\Gr_F^p\ne 0$).
Corollary~5 also implies that $N^{n-1}\ne 0$ on $H^{n-1}(X_{\infty},
\Q)$ for $r\ge n+1$, since
$$N^k:\Gr^W_{n-1+k}H^{n-1}(X_{\infty},\Q)\simto
\Gr^W_{n-1-k}H^{n-1}(X_{\infty},\Q)(-k)\q(k>0).$$
For $j\ne n-1$, we have $H^j(X_{\infty},\Q)\cong\Q(-j/2)$ or $0$
in the case of Corollary~5.
It is easy to show that the direct sums in Corollary~5 are direct
factors of $P_N\Gr^W_{n-1+k}H^{n-1}_{\prim}(X_{\infty},\Q)$.
Using the $N$-primitive decomposition, Corollary~5 is then
equivalent to the numerical equality
$$\h{$P_n(d_0)=\sum_{|I|\ge 1}|I|\binom{r}{|I|}P_n(\d_I)+
\sum_{|I|\ge 2}(|I|-1)\binom{r}{|I|}P_n(\d'_I)+
n\binom{r-1}{n}$},\leqno(0.4)$$
where $\d_I=(d_k)_{k\in I}$, $\d'_I=(d_k)_{k\in I\cup\{0\}}$, and
$P_n(\d)$ denotes the dimension of the reduced middle primitive
cohomology of a smooth complete intersection of multidegree
$\d\in\Z_{>0}^k$ in $\P^n$.
It is possible to show directly (0.4) in some simple cases,
see Remark~(3.5)(i).
By Corollary~5, we can calculate the limit mixed Hodge structure
without using a blow-up of $X$ in this case.
It generalizes a calculation of the limit mixed Hodge structure
for $n=4$, $d_0=r=5$ in [11] where a blow-up of $X$ is used.
However, we cannot calculate the group of connected components of
a fiber of the N\'eron model as in loc.~cit.\ by using our method.

As an application of Theorem~5 where $r=2$, we get the following

\ms\nin
{\bf Theorem~6.} {\it
Let $Y$ be a smooth complex projective variety of dimension
$n\ge 2$.
Let $L$ be a very ample line bundle defining a closed embedding
$Y\into\P^N$.
Assume $L$ is a tensor product of $k$ very ample line bundles
where $k=3$ if $\,Y$ is projective space of even dimension,
and $k=2$ otherwise.
Then the vanishing cycles of a Lefschetz pencil do not vanish,
i.e.\ the restriction morphism $i_s^*:H^{n-1}(Y,\Q)\to
H^{n-1}(Y_s,\Q)$ is non-surjective where $i_s:Y_s\to Y$ is a
general member of the linear system $|L|=(\P^N)^{\vee}$.}

\ms
This is an improvement of [15], Cor.~6.4 where the assertion was
shown for $L=L'^{\otimes d}$ with $L'$ ample and $d\gg 1$.
It can be used to show the vanishing of some direct factor of
the decomposition in [4] for the direct image of the constant
sheaf by $\Y\to(\P^N)^{\vee}$ where $\Y$ is the total space
of the universal family of the hyperplane sections of
$Y\subset\P^N$.
Note that non-surjectivity of $i_s^*$ implies that the
discriminant of $\pi$ (i.e.\ the dual variety of $Y$ in
$(\P^N)^{\vee}$) has codimension 1.
The converse is true in the $n$ odd case,
using the fact that the eigenvalue of the local monodromy of the
vanishing cycle is $-1$ by the Picard-Lefschetz formula [17],
see [15], Th.~6.3.
Note also that taking the tensor product of two very ample line
bundles corresponds to the composition with the Segre embedding.
It is possible to prove Theorem~6 by using some arguments
in [13], [20].

Recently we are informed that our paper is closely related to some
results in [1] and [2].

\ms
In Section~1, we recall some basics of mixed Hodge modules and
show a lemma used in the proof of Theorem~2.
In Section~2, we prove Theorems~1--4 and Corollaries~3--4.
In Section~3, we introduce quasi-semistable degenerations and
prove a generalization of Theorem~5 and also Corollary~5.
In Section~4, we show Theorem~6.

\section{Preliminaries}

\nin
{\bf 1.1.~Weight filtration.}
Every mixed Hodge module $M$ on a complex analytic space $X$
has a canonical weight filtration $W$ in the category of mixed
Hodge modules $\MHM(X)$, and every morphism of mixed Hodge modules
is strictly compatible with the weight filtration $W$.
We say that $M$ is pure of weight $n$ if $\Gr^W_kM=0$ for $k\ne n$.
If a mixed Hodge module $M$ is pure, then it has a strict support
decomposition $M=\mopl_Z\,M_Z$,
where $Z$ runs over the irreducible closed analytic subspaces of
$X$, and $M_Z$ has strict support $Z$, i.e. its support is $Z$
and there is no nontrivial sub nor quotient object with strictly
smaller support.
Pure Hodge modules are semisimple since they are assumed to be
polarizable, see [22], 5.1--2.

\ms\nin
{\bf 1.2.~Nearby and vanishing cycle functors.}
With the above notation, assume $X$ is pure dimensional, and
let $X_i$ be the irreducible components of $X$.
Let $M$ be a pure Hodge module of weight
$n$ which is a direct sum of pure Hodge modules $M_{X_i}$ with
strict support $X_i$.
Let $f$ be a holomorphic function on $X$ which is nonconstant
on any $X_i$.
By definition ([22], 5.1.6) the weight filtration $W$ on the nearby
and vanishing cycles $\psi_fM,\ph_{f,1}M$ is the monodromy
filtration shifted by $n-1$ and $n$ respectively.
So we have
$$\aligned
N^k:\Gr^W_{n-1+k}\psi_fM&\simto(\Gr^W_{n-1-k}\psi_fM)(-k),\\
N^k:\Gr^W_{n+k}\ph_{f,1}M&\simto(\Gr^W_{n-k}\ph_{f,1}M)(-k).\endaligned
\leqno(1.2.1)$$
As for the non-unipotent monodromy part, we have
$$(\psi_{f,\ne1}M,W)=(\ph_{f,\ne1}M,W).$$

Set $Y=f^{-1}(0)$ with the inclusion $i:Y\into X$.
Since $M$ has strict support $X$ and $Y\ne X$, we have
$$\HH^ji^*M=0\q\h{for}\q j\ne-1,$$
and there is a short exact sequence of mixed Hodge modules on $Y$
$$0\to\HH^{-1}i^*M\to\psi_{f,1}M\to\ph_{f,1}M\to 0.\leqno(1.2.2)$$
In fact, the functor $i^*$ is defined by the mapping cone of
$\psi_{f,1}\to\ph_{f,1}$, i.e. we have a distinguished triangle
$$i^*[-1]\to\psi_{f,1}\to\ph_{f,1}\to.\leqno(1.2.3)$$
This is slightly different from the usual one since
$\psi_{f,1},\ph_{f,1}$ in this paper are shifted by $-1$ so that
they preserve perverse sheaves.
We have moreover isomorphisms
$$\HH^{-1}i^*M=\Ker\,N\subset\psi_{f,1}M,\q
\ph_{f,1}M=\Coim\,N,\leqno(1.2.4)$$
where $N:\psi_{f,1}M\to\psi_{f,1}M(-1)$, see [22], 5.1.4.
Indeed, we have a surjection and an injection
${\rm can}:\psi_{f,1}M\onto\ph_{f,1}M$ and
${\rm Var}:\ph_{f,1}M\into\psi_{f,1}M(-1)$
such that $N={\rm Var}\scirc{\rm can}$.

\ms\nin
{\bf 1.3.~Constant sheaf case.}
If $X$ is an intersection of hypersurfaces in a complex manifold
$V$ (shrinking $X$ if necessary), then we have $\Q_X$ in the
derived category of mixed Hodge modules $\DD:=D^b\MHM(X)$ as in
the proof of Prop.~2.19 in [23].
Indeed, if $i_X:X\into V$ and $i_j:g^{-1}(0)\to V$ denote
the inclusions where $\bigcap_jg^{-1}(0)=X$, then
$$(i_X)_*i_X^*=\mpro_j(i_j)_*i_j^*\q\h{with}\q
(i_j)_*i_j^*=C(\psi_{g_j,1}\to\ph_{g_j,1}).$$

Note that (0.1) is equivalent to
$$\aligned\Hom_{\DD}(\Gr^W_n\HH^n\Q_X,M)&=\Hom_{\DD}(\HH^n\Q_X,M)
=\Hom_{\DD}(\Q_X[n],M)\\
&=\Ext^{\dim Z-n}_{\DD}(\Q_Z[\dim Z],M)=0\endaligned$$
for any pure Hodge module $M$ of weight $n$ with
$Z:=\supp M\subset\Sing X$
(where we may assume that $Z$ is an intersection of hypersurfaces
shrinking $V$ if necessary).
Indeed, the strict support decomposition implies
$$\Gr^W_n\HH^n\Q_X=\IC_X\Q\oplus M',$$
where $M'$ is a pure Hodge module of weight $n$ and
$\supp M'\subset\Sing X$,
and it is enough to show that $M'=0$.

For a holomorphic function $f$ on $X$ and $x\in Y:=f^{-1}(0)$,
let $i_x:\{x\}\into Y$ denote the inclusion, and $F_x$ denote
the Milnor fiber around $x$.
Since $\psi_f$ and $\ph_f$ in this paper are shifted by $-1$,
we have
$$\aligned H^ji_x^*\psi_f(\Q_X[n])&=H^{n-1+j}(F_x,\Q),\\
H^ji_x^*\ph_f(\Q_X[n])&=\tH^{n-1+j}(F_x,\Q).\endaligned
\leqno(1.3.1)$$

The following will be used in the proof of Theorem~2.

\ms\nin
{\bf Lemma~1.4.} {\it Let $\A$ be the category consisting of
$(M_{\ssb},N)$ where $M_{\ssb}=\mopl_{k\in\Z}M_k$ with
$M_k$ a pure Hodge module of weight $k$ and
$N:M_{\ssb}\to M_{\ssb}(-1)$ is a morphism of graded Hodge modules
\(here $(M_{\ssb}(-1))_k:=M_{k-2}(-1).)$
Morphisms of $\A$ are morphisms of graded Hodge modules compatible
with the action of $N$.
Assume there is a commutative diagram of exact sequences in $\A$
$$\begin{matrix}&&0&&0&&C'_{\ssb}\\
&&\da&&\da&&\bigcap\\
0&\to&A'_{\ssb}&\to&B_{\ssb}&\to&C_{\ssb}&\to&0\\
&&\bigcap&&||&&\don\\
0&\to&A_{\ssb}&\to&B_{\ssb}&\to&C''_{\ssb}&\to&0\\
&&\don&&\da&&\da\\
&&A''_{\ssb}&&0&&0\end{matrix}$$
such that
$$\aligned&N^k:B_{n-1+k}\simto B_{n-1-k}(-k)\q(k>0),\\
&A_{\ssb}=\Ker(N:B_{\ssb}\to B_{\ssb}(-1)).\endaligned$$
Then, choosing a graded splitting of $A_{\ssb}\onto A''_{\ssb}$,
we have an isomorphism in $\A$
$$C_{\ssb}=(\mopl_{k\ge 1}\mopl_{j=1}^k\,A'_{n-1-k}(-j))\oplus
(\mopl_{k\ge 0}\mopl_{j=0}^k\,A''_{n-1-k}(-j)),$$
such that $N:C_{\ssb}\to C_{\ssb}(-1)$ is identified with a
morphism induced by the identity on $A'_{n-1-k}(-j)$,
$A''_{n-1-k}(-j)$.}

\ms\nin
{\it Proof.}
Note first that $C'_{\ssb}=A''_{\ssb}$ by the snake lemma.
Since $A_{\ssb}$ is identified with the $N$-primitive part of
$B_{\ssb}$ up to Tate twists, we have the Lefschetz decompositions
$$B_{\ssb}=\mopl_{k\ge 0}\mopl_{j=0}^k\,A_{n-1-k}(-j),\q
C''_{\ssb}=\mopl_{k\ge 1}\mopl_{j=1}^k\,A_{n-1-k}(-j),$$
such that $N:B_{\ssb}\to B_{\ssb}(-1)$ is identified with a
morphism induced by the identity on $A_{n-1-k}(-j)$,
and similarly for $C''_{\ssb}$.
Then, choosing a graded splitting of $A_{\ssb}\onto A''_{\ssb}$,
the assertion follows.

\ms\nin
{\bf Remark~1.5.}
It is well-known to the specialists that the condition
$\Q_X[n]=\IC_X\Q$ is equivalent to that $X$ is a $\Q$-homology
manifold.
Indeed, the former condition implies that
$$\Q_X[n]=(\D\Q_X)(-n)[-n],\leqno(1.5.1)$$
using the self-duality $\D(\IC_X\Q)=\IC_X\Q(n)$ where $\D$
denotes the functor associating the dual.
Then (1.5.1) implies that $X$ is a $\Q$-homology manifold.
The converse is easy, see also [10], [18].

\section{Proof of Theorems 1--4 and Corollaries 3--4}

\nin
{\bf 2.1.~Proof of Theorems~1 and 4.}
Set $\G_X=C(\Q_X\to(\IC_X\Q)[-n])[-1]$ so that we have a
distinguished triangle
$$\G_X\to\Q_X\to(\IC_X\Q)[-n]\to.\leqno(2.1.1)$$
Since the $\HH^j\G_X$ are supported on $Y$, $\G_X$ can be
identified with a complex of mixed Hodge modules on $Y$ (i.e.\
it is viewed as an abbreviation of $i_*\G_X$), see [23], 2.23.
Applying the functor $i^*$, we get then a distinguished
triangle on $Y$
$$\G_X\to\Q_Y\to i^*(\IC_X\Q)[-n]\to.\leqno(2.1.2)$$
Since $f$ is nonconstant on any $n$-dimensional irreducible
component of $X$, we have
$$\HH^ji^*(\IC_X\Q)=0\q(j\ne -1).\leqno(2.1.3)$$
By (2.1.1) and (0.1--2) we have
$$\HH^j\G_X\simto\HH^j\Q_X\,\,(j<n),\,\,\,\,
\HH^n\G_X\simto W_{n-1}\HH^n\Q_X,\,\,\,\,\HH^j\G_X=0\,\,(j>n),
\leqno(2.1.4)$$
and Theorem~1 follows by using the long exact sequence associated
to the distinguished triangle (2.1.2).
Then Theorem~4 follows from (1.2.4) and Theorem~1 together with
the semisimplicity of pure Hodge modules.

\ms\nin
{\bf 2.2.~Proof of Theorem~2.}
Applying (1.2.3) to (2.1.1) shifted by $n$, and taking the
associated long exact sequences, we get a commutative diagram of
exact sequences
$$\begin{matrix}0&&0&&&&0\\\da&&\da&&&&\da\\
\HH^{n-1}\G_X&\simto&\HH^{n-1}\G_X&\to&0&\to&\HH^n\G_X\\
||&&\da&&\da&&\da\\
\ph_{f,1}\HH^{n-1}\Q_X&\into&\HH^{n-1}\Q_Y&\to&\psi_{f,1}\HH^n\Q_X&
\onto&\ph_{f,1}\HH^n\Q_X\\\da&&\da&&||&&\da\\
0&\to&\HH^{-1}i^*(\IC_X\Q)&\into&\psi_{f,1}(\IC_X\Q)&\onto&
\ph_{f,1}(\IC_X\Q)\\\da&&\don&&\da&&\da\\
\HH^n\G_X&\simto&\HH^n\G_X&\to&0&&0\end{matrix}$$
where all the squares commute since we have the vanishing of
certain terms of the the squares of anti-commutativity (see
[3], 1.1.11).
Here $\G_X$ is identified with $i_*\G_X$ so that
$\ph_{f,1}\G_X$ and $i^*\G_X$ are identified with $\G_X$.
The fourth row for $\IC_X\Q$ is a short exact sequence since
$f$ is nonconstant on any irreducible component $X_i$ of $X$ and
$\IC_{X_i}\Q$ has strict support $X_i$, see (1.2.2).
By (2.1.4) the above diagram induces a diagram of the snake lemma
$$\begin{matrix}&&0&&0&&W_{n-1}\HH^n\Q_X\\&&\da&&\da&&\bigcap\\
0&\to&\HH^{n-1}\Q_Y/\HH^{n-1}\Q_X&\to&\psi_{f,1}\HH^n\Q_X&\to&
\ph_{f,1}\HH^n\Q_X&\to&0\\&&\bigcap&&||&&\don\\
0&\to&\HH^{-1}i^*(\IC_X\Q)&\to&\psi_{f,1}(\IC_X\Q)&\to&
\ph_{f,1}(\IC_X\Q)&\to&0\\&&\don&&\da&&\da\\
&&W_{n-1}\HH^n\Q_X&&0&&0\end{matrix}$$
The weight filtration $W$ on $\psi_{f,1}(\IC_X\Q)$ and
$\ph_{f,1}(\IC_X\Q)$ are the monodromy filtrations shifted by
$n-1$ and $n$ respectively since $\IC_X\Q$ is pure of weight $n$,
see (1.2.1).
Moreover, we have by (1.2.4)
$$\HH^{-1}i^*(\IC_X\Q)=\Ker\,N\subset\psi_{f,1}(\IC_X\Q).$$
So the assertion follows from Lemma~(1.4).

\ms\nin
{\bf 2.3.~Proof of Theorem~3.}
Since $\tH^{n-1+j}(F_x,\Q)_1=\H^j(\ph_{f,1}\Q_X[n])_x$,
Theorem~3 follows from Theorem~2 together with (0.3) for $X,Y$.

\ms\nin
{\bf 2.4.~Proof of Corollaries~3 and 4.}
The graded pieces $$V_k:=\Gr^W_k\ph_{f,1}\1\HH^0(\Q_X[n])$$ have
strict support decompositions
$$V_k=(\mopl_i(V_k)_{Z_i})\oplus(\mopl_{z\in Z}(V_k)_{\{z\}}),$$
where $Z_i$ are the irreducible components of $Z$.
By the assumption on the local monodromy $T_{i,z}$ along $Z_i$
around $z$, we have
$$H^ji_z^*(V_k)_{Z_i}=0\q\h{for any}\,\,\,j,$$
since it is calculated by the mapping cone of the action of
$T_{i,z}-id$.
On the other hand we have clearly
$$H^ji_z^*(V_k)_{\{z\}}=0\q\h{for}\,\,\,j\ne 0,$$
These imply the $E_1$-degeneration of the spectral sequence
associated to the functor $H^{\ssb}i_z^*$ and the filtration $W$
on $\HH^0\ph_{f,1}(\Q_X[n])$, i.e.
$$\h{$H^{\ssb}i_z^*$ and $\Gr^W_{\ssb}$ commute on
$\HH^0\ph_{f,1}(\Q_X[n])$.}$$
We have moreover
$$H^ji_z^*(\HH^0\ph_{f,1}(\Q_X[n]))=\begin{cases}H^{n-1}(F_z,\Q)_1&
\text{if $\,j=0$,}\\ 0&\text{if $\,j\ne 0$.}\end{cases}$$
Indeed, the assertion for $i=0$ follows from the spectral sequence
$$\aligned E_2^{i,j}=H^ii_z^*(\HH^j\ph_{f,1}(\Q_X[n]))\Rightarrow
H^{i+j}i_z^*(\ph_{f,1}(\Q_X[n]))\\ =H^{i+j+n-1}(F_z,\Q)_1,
\endaligned$$
since $E_2^{i,j}=0$ if $i>0$ or $i<0,j=0$.
We have a similar spectral sequence with $\ph_{f,1}(\Q_X[n])$
replaced by $\G_X$, since $\Gr^W_k\HH^n\G_X$ is a direct factor
of $V''_k$ by Theorem~2.
So we get similarly
$$H^ji_z^*(\HH^n\G_X)=\begin{cases}H^ni_z^*\G_X&
\text{if $\,j=0$,}\\ 0&\text{if $\,j\ne 0$.}\end{cases}$$
Here we have
$$H^ni_z^*\G_X=H^{-1}i_z^*\IC_X\Q,$$
since $n\ge 2$ (because $\dim Z=1$).
In case $\HH^{n-1}\Q_X=0$, we have a similar argument with
$X,n$ replaced by $Y,n-1$, since we have the isomorphism
$$K'_k=(\Gr^W_{k-2}\HH^{n-1}\Q_Y)(-1)\,\,\,(k\le n).$$
So Corollaries~3--4 follow from Theorem~2.

\ms\nin
{\bf Remarks~2.5.}
(i) In the case $X$ is a $\Q$-homology manifold, Theorem~3
implies that $Y$ is a $\Q$-homology manifold (i.e.\
$\Q_Y[n-1]=\IC_Y\Q$) if and only if $\ph_{f,1}\Q_X=0$.
This seems to be known to the specialists at least if $X$ is
smooth, see also [28].
Note that this equivalence does not hold for the singular case
even if $X$ is a complete intersection so that
$\HH^i(\Q_X[n])=\HH^i(\Q_Y[n-1])=0$ for $i<0$, see Theorem~3.

\ms
(ii) If $x$ is an isolated singularity of $X,Y$, let
$L_{X,x}$ denote the link of $(X,x)$, and similarly for $Y$.
In this case $\G_X$ in (2.1.1) is supported on $\{x\}$, and we
have isomorphisms for $j\le n$
$$\HH^j\G_X=\H^j\G_X=H^{j-1}((\IC_X\Q[-n])_x/\Q)=
\tH^{j-1}(L_{X,x},\Q)=H^{j}_{\{x\}}\Q_X.$$
Combining this with (2.1.4), we see that $\HH^j\Q_X$ and
$\HH^j\Q_Y$ in Theorem~1 can be replaced with the local
cohomology groups in this case, see also [21], [27].

\ms
(iii) If $X$ is a complete intersection and $x$ is an isolated
singularity of $X,Y$, then it is also possible to prove
Theorem~3(c) as follows.
We have an exact sequence with the notation of Remark~(ii) above
$$0\to H^{n-2}(L_{Y.x})\to H_c^{n-1}(L_{X.x}\setminus L_{Y,x})\to
H^{n-1}(L_{X.x})\buildrel{i'^*}\over\to H^{n-1}(L_{Y.x}),$$
where $i':L_{Y.x}\to L_{X,x}$ is the inclusion,
and the morphism $i'^*$ vanishes since $H^{n-1}(L_{X.x})$ has
weights $\le n-1$ and $H^{n-1}(L_{Y.x})$ has weights $>n-1$,
see e.g.\ [9].
(It does not seem easy to prove this vanishing without using Hodge
or $\ell$-adic theory.)
Then the assertion follows from the Wang sequence associated to
the Milnor fibration $L_{X.x}\setminus L_{Y,x}\to S^1$ constructed
in [14].

\section{Case of quasi-semistable degenerations}

\nin
In this section we introduce quasi-semistable degenerations, and
prove a generalization of Theorems~5 and also Corollary~5.

\ms\nin
{\bf 3.1.~Quasi-semistable degenerations.}
Let $f:X\to S$ be a proper morphism of complex analytic spaces
such that $S$ is an open disc.
We say that $f$ is a quasi-semistable degeneration if
$X_0:=f^{-1}(0)$ is a reduced variety with simple normal crossings,
$(X,x)$ for each $x\in X_0\cap\Sing X$ is analytically isomorphic to
$$(h^{-1}(0),0)\subset(\C^{n+1},0)\q\h{with}\q
h=y_1\cdots y_k-y_nt,\leqno(3.1.1)$$
and moreover $f$ is locally identified with $t$ by choosing a local
coordinate of $S$.
Here $y_1,\dots,y_n,t$ are the coordinates of $\C^{n+1}$, and
$k\in[2,n-1]$ may depend on $x$.
Note that we have on a neighborhood of $x\in X_0\cap\Sing X$
$$\Sing X=\Sing X_0\cap\{y_n=0\}\subset X_0.\leqno(3.1.2)$$
Let $Y_k\,(1\le k\le r)$ be the irreducible components of $X_0$
which are smooth by hypothesis.
We assume $X_0$ algebraic (or K\"ahler in a generalized sense
that each $Y_i$ has a K\"ahler form such that their cohomology
classes come from a cohomology class on $X_0$). Set
$$Y_I=\mcap_{k\in I}Y_k,\q Y'_I=Y_I\cap\Sing X\q\h{for}\q
I\subset\{1,\dots,r\}.$$
When
we consider $Y'_I$, we will assume $|I|\ge 2$ so that
$Y_I\subset\Sing X_0$ and
$$Y'_I=Y_I\cap\{y_n=0\}\q\h{locally}.$$
Since $X$ is locally a hypersurface, $\Q_X[n],\Q_{X_0}[n-1]$ are
perverse sheaves.
By (3.1.2) we have
$$\psi_{f,1}\IC_X\Q=\psi_{f,1}(\Q_X[n])=\psi_f(\Q_X[n]),
\leqno(3.1.3)$$
where the last isomorphism follows from the fact that the
monodromy on the Milnor cohomology of $f$ at $x\in X_0$ is the
identity.
Indeed, with the notation of (3.1.1), we have a
geometric monodromy induced by the action of $\R$ defined by
$$\alpha:(y_1,\dots,y_n,t)\mapsto(e^{2\pi i\alpha}y_1,y_2,
\dots,y_n,e^{2\pi i\alpha}t),$$
and it is the identity for $\alpha=1$.

The advantage of quasi-semistable degeneration is that it is
quite easy to construct examples as is shown in Theorem~5 and
Corollary~5.
To get further a (nonreduced) semistable model we would need
blowing-ups as in [11], [13].

\ms\nin
{\bf Theorem~3.2.} {\it The conclusion of Theorem~$5$ holds for
any quasi-semistable degeneration.}

\ms\nin
{\it Proof.}
By (3.1.3) the spectral sequence is associated with the weight
filtration $W$ on $\psi_{f,1}\IC_X\Q$ which is the monodromy
filtration shifted by $n-1$.
By [23], 2.14 (together with the definition of $\psi_t$),
we have isomorphisms as mixed Hodge structures
$$H^{j}(X_0,\psi_f(\Q_X[n]))=
\psi_t{}^p\!R^jf_*(\Q_X[n]))=H^{j+n-1}(X_{\infty},\C),$$
where $\psi_f,\psi_t$ are shifted by $-1$ so that they
preserve perverse sheaves and mixed Hodge modules.
By Theorem~4, it is then sufficient to show the following
isomorphisms for $k<n$
$$\aligned\Gr^W_k(\Q_{X_0}[n-1])&=\mopl_{|I|=n-k}\,\Q_{Y_I}
[\dim Y_I],\\\Gr^W_k(\Q_X[n])&=\mopl_{|I|=n-k+1}\,\Q_{Y'_I}
[\dim Y'_I](-1).\endaligned\leqno(3.2.1)$$
The first isomorphisms are well-known.
For the second isomorphisms, we use the function $h$ in (3.1.1)
defining locally $X$.
By (1.2.4) and using the $N$-primitive decomposition, there are
isomorphisms for $k<n$
$$\aligned\Gr^W_k(\Q_X[n])&=\Gr^W_k\psi_{h,1}(\Q_{X'}[n+1])\cap
\Ker\,N\\&=\Gr^W_{k+2}\ph_{h,1}(\Q_{X'}[n+1])(1)\cap\Ker\,N.
\endaligned$$
The assertion is then reduced to the first isomorphisms of (3.2.1)
by using Lemma~(3.3) below.
More precisely, it proves the assertion locally in classical
topology, and we have to show that the direct factors of
$\Gr^W_k(\Q_X[n])\,(k<n)$ are {\it globally} constant sheaves.
This is inductively reduced to the case $k=n-1$.
(Indeed, the local extension class between the $\Gr^W_k(\Q_X[n])$
induces the restriction morphisms $\Q_{Y'_I}\to\Q_{Y'_J}$ for
$I\subset J$ with $|I|=|J|-1$, see also Remark~(3.5)(iv) below.)

In case $k=n-1$, we may assume $r=2$ replacing $Y$ with
an affine open subvariety if necessary.
In this case, $X$ is equisingular along $Y'_{\{1,2\}}$ and has an
ordinary double point by restricting to a transversal slice to
$Y'_{\{1,2\}}$.
So it is enough to show that the restriction of $\H^j\IC_X$ to
$Y'_{\{1,2\}}$ is a constant sheaf for any $j$.
Then the assertion is proved by taking the blow-up of $X$ along
$Y'_{\{1,2\}}$ which gives a resolution of singularities
and using the direct image of the constant sheaf since the
intersection complex is a direct factor of the direct image
by the decomposition theorem~[3].

\ms\nin
{\bf Lemma~3.3.} {\it
Let $Z$ be a complex manifold, and $g$ be a holomorphic function
on $Z$ such that $g^{-1}(0)$ is a reduced divisor with simple
normal crossings.
Let $X'=Z\times\C^2$, and $h=g+z_1z_2$ on $X'$ where $z_1,z_2$
are the coordinates of $\C^2$.
Then we have a canonical isomorphism compatible with the action
of $N$}
$$\ph_{h,1}(\Q_{X'}[\dim X'])=\ph_{g,1}(\Q_Z[\dim Z])(-1).$$

\ms\nin
{\it Proof.}
This is a special case of the Thom-Sebastiani theorem for
Hodge modules which was shown in an unpublished manuscript of
the second author which was typed at RIMS in 1990.
For the proof of Theorem~5 it is enough to show it for the
underlying perverse sheaves with $\C$-coefficients since the
$N$-primitive part of each graded piece is a direct sum of the
constant sheaves on the strata in this case.
Then the assertion also follows from [24], 4.1.

\ms\nin
{\bf 3.4.~Proof of Corollary~5.}
By Theorem~5 the $E_1$-term $E_1^{-i,j+i}$ is the direct
sum of
$$H^{j-i+n-1-2l}(Y_I)(-i-l)\q\h{and}\q H^{j-i+n-1-2l-2}(Y'_{I'})
(-i-l-1)$$ over $I$, $I'$ satisfying respectively the conditions
$$|I|=i+2l+1,\q |I'|=i+2l+2,$$
where $l$ satisfies the conditions $i+l\ge 0$, $i\ge 0$, and
$$\begin{matrix}j\le n-1-i-2l,& j\ge -n+1+i+2l & \h{for}&\!\!\!Y_I,
&\h{$\qquad\quad\,$}\\
j\le n-1-i-2l-2,& j\ge -n+1+i+2l+2 & \h{for}&\!\!\!
{}^{\raise8pt\h{$\,$}}Y'_{I'}.\\ \end{matrix}
\leqno(3.4.1)$$
Note that $\dim Y_I=n-1-i-2l$ and
$\dim Y_{I'}=n-1-i-2l-2$.

Let $'\!H^{\ssb}(Y_I)$ denote the orthogonal complement of
the middle primitive cohomology in $H^{\ssb}(Y_I)$ if
$\dim Y_I\ne 0$, and the image of the canonical morphism
$\Q\to H^0(Y_I)$ if $\dim Y_I=0$.
Note that $'\!H^j(Y_I)=H^j(Y_I)$ if $j\ne\dim Y_I$.
Define $'\!H^{\ssb}(Y'_I)$ similarly.
The differential $d_1^{p,q}$ is expressed by using the restriction
and Gysin morphisms up to constant multiples.
(This is shown by using the extension classes at the level of
sheaves, see Remark~(3.5)(iv) below and also [20].)
Then $'\!H^{\ssb}(Y_I)$ and $'\!H^{\ssb}(Y'_I)$ define a
subcomplex $\E_1^{\ssb,\ssb}$ of $E_1^{\ssb,\ssb}$.
So it is enough to show that its cohomology $\E_2^{-i,j+i}$
satisfies
$$\E_2^{-i,i}={\buildrel{m}\over\mopl}\,\Q((1-n-i)/2)\q\h{if}
\,\,\,i+n-1\,\,\h{is even and}\,\,\,|i|\le n-1,\leqno(3.4.2)$$
where $m=\binom{r-1}{n}$.
Note that $\E_2^{-i,i+j}=0$ for $i+j+n-1$ odd.
(Using (3.4.2), the assertion can be reduced to the case $d_k=1$
for any $k>0$.
Indeed, $\E_1^{\ssb,\ssb}$ is independent of the $d_k$ as long as
$r$ is fixed, where the differential can be neglected by using
the Euler characteristic as below.)

Take $k\in\Z$, and consider the subcomplex $\E_1^{\ssb-k,k}$.
This is defined by the direct factors of $\E_1^{-i,i+j}$ with
$$i+j=k.$$
Here we may assume $k-n-1$ is even since the complex vanishes
otherwise.
It is enough to calculate the Euler characteristic of
$\E_1^{\ssb-k,k}$ since we have for $j\ne n-1$
$$\h{$H^j(X_{\infty},\Q)=\Q(-j/2)\,\,$ if $\,j\,$ is even, and
$\,0\,$ otherwise.}$$
By the self-duality of $(\psi_{f,1},N)$, the $E_1$-complex is
self-dual, and we may assume $k\ge 0$.

Let $I^l$ denote the image filtration on $\E_1$
defined by $\Im\,N^l$.
The action of $N$ on $\E_1$ is defined so that
the index $l$ increases by 1 and $\E_1/I^1$ is isomorphic to the
primitive part.
Thus the index $l$ is constant on the graded complex
$\Gr_I^l\,\E_1^{\ssb-k,k}$.
We will use the index $i$ instead of $j$ so that the complex
will be denoted by $\Gr_I^l\,\E_1^{-\ssb,k}$
(because of the relation $|I|=i+2l+1$, etc.)
Note that the differential decreases the index $i$ by 1.

Since it is sufficient to calculate the Euler characteristic of
$\Gr_I^l\,\E_1^{-\ssb,k}$, we may modify the differential as we
like.
So we may calculate it separately for $Y_I$ and $Y'_I$.
We first consider the complex consisting only of the cohomology
of $Y_I$, which will be denoted by $(\Gr_I^l\,\E_1^{-\ssb,k})_Y$.
Let $K^{(r)}_{\ssb}$ denote the Koszul complex associated to $r$
morphisms which are the identity on $\Q$, where
$\dim K^{(r)}_i=\binom{r}{i}$.
Since $|I|=i+2l+1$ corresponds to the index $i$ of the Koszul
complex up to a shift, we may assume (by modifying the differential
of $(\Gr_I^l\,\E_1^{-\ssb,k})_Y$ if necessary)
$$(\Gr_I^l\,\E_1^{-\ssb,k})_Y\cong(\sigma_{\le q}\sigma_{\ge p}
K^{(r)}_{\ssb})[-2l-1],$$
where $\sigma_{\ge p}$ is the truncation which preserves the
the components of degree $\ge p$ and replace it with 0 for
the degree $<p$ (and similarly for $\sigma_{\le q}$), see [5].
Note that for $p<q$
$$\h{$\dim H_i(\sigma_{\le q}\sigma_{\ge p}K^{(r)}_{\ssb})=
\begin{cases}\binom{r-1}{p-1}&\text{if $i=p$,}\\
\binom{r-1}{q}&\text{if $i=q^{\raise10pt\hbox{$\,$}}$,}\\
\,\,\,\,\,\,0^{\raise8pt\hbox{$\,$}}&\text{otherwise.}\end{cases}$}
\leqno(3.4.3)$$
(This follows from the well-known formula
$\binom{r}{i}=\binom{r-1}{i}+\binom{r-1}{i-1}$ together with
the acyclicity of $K^{(r)}_{\ssb}$ by decomposing it into the
short exact sequences using $\Ker\,d=\Im\,d$.)

The numbers $p,q$ are determined by using the conditions in
(3.4.1) for $Y_I$, and we get
$$p=l+1,\q q=l+(n+1+k)/2.$$
Moreover the range of $l$ is given by
$$0\le l\le A:=(n-1-k)/2,$$
which implies $k\le n-1$.
Then we get the ranges of $p,q$ (when $l$ varies)
$$1\le p\le A+1,\q A+k+1\le q\le n.$$

We apply a similar argument to $(\Gr_I^l\,\E_1^{-\ssb,k})_{Y'}$
which consists of the cohomology of $Y'_I$.
Since $|I'|=i+2l+2$, we have
$$(\Gr_I^l\,\E_1^{-\ssb,k})_{Y'}\cong(\sigma_{\le q}\sigma_{\ge p}
K^{(r)}_{\ssb})[-2l-2],$$
where the shift of the Koszul complex is different from the above
complex by 1.
Applying further a similar argument we get
$$\aligned p=l+2,\q q=l+(n+1+k)/2,\q 0\le l\le A-1,\\
2\le p\le A+1,\q A+k+1\le q\le n-1,\q k\le n-3.\endaligned$$
If $k\le n-3$, we take the Euler characteristic of
$\E_1^{-\ssb,k}$, and get cancelations for
$$2\le p\le A+1\q\h{and}\q A+k+1\le q\le n-1.$$
Thus, using (3.4.3), only $\binom{r-1}{0}$ for $p=1$ and
$\binom{r-1}{n}$ for $q=n$ remain .
The first term corresponds to $H^{n-1+k}(X_{\infty})$ if $k>0$,
and to the orthogonal complement of the middle primitive part
if $k=0$.
So we get the desired assertion.
If $k=n-1$, then the range of $l$ consists only of $\{0\}$ for
$(\Gr_I^l\,\E_1^{-\ssb,k})_Y$, and is empty for
$(\Gr_I^l\,\E_1^{-\ssb,k})_{Y'}$.
So we get the same conclusion.

\ms\nin
{\bf Remarks~3.5.} (i)
For the proof of Corollary~5 we may assume $d_k=1$ for any $k>0$,
see a remark after (3.4.2).
Set $d:=d_0=r$.
Let $e_n(d)$ denote the Euler number of a smooth hypersurface $Z$
of degree $d$ in $\P^n$ for $n,d\ge 2$ (e.g.\ $Z=Y_0$ or $Y'_I$).
It is well-known that
$$\h{$e_n(d)=-\sum_{i=0}^{n-1}\binom{n+1}{i}(-d)^{n-i}=
n+1+((1-d)^{n+1}-1)/d$}.$$
This can be verified by using an inductive formula
$e_n(d)=nd-(d-1)e_{n-1}(d)$ which is easily shown for Fermat
hypersurfaces.
For $n\ge 1$, set
$$P_n(d):=\dim\tH^{n-1}_{\prim}(Z)=(d-1)((d-1)^n-(-1)^n)/d,$$
where $H^{n-1}_{\prim}(Z)$ denotes the primitive cohomology, and
the last equality for $n\ge 2$ follows from the above formula.
By the argument using $\E_1^{\ssb,\ssb}$ in (3.4), the direct
sums in Corollary~5 are direct factors of the $N$-primitive
part.
So Corollary~5 is equivalent to the numerical equality (0.4).
In the case $d_k=1\,(k>0)$ and $d_0=d$, the latter becomes
$$\h{$P_n(d)=\sum_{0\le k\le n-3}(k+1)\binom{d}{k+2}
P_{n-k-2}(d)+n\binom{d-1}{n}$}.$$
It is possible to prove this by a direct calculation if $n$ is
quite small (e.g. if $n\le 5$).

\ms
(ii) Let $C(n+1,d,j)$ denote the integers such that
$$(t+\cdots+t^{d-1})^{n+1}=\msum_{j=n+1}^{(n+1)(d-1)}\,
C(n+1,d,j)t^j.$$
Up to the multiplication by $t^{n+1}$, this coincides with the
Poincare polynomial of the graded module
$$\C[x_0,\dots,x_n]/(x_0^{d-1},\dots,x_n^{d-1}).$$
Using a Koszul complex which is associated with the multiplication
by $x_i^{d-1}$ and gives a free resolution of the above module,
we get as is well-known to the specialists
$$\h{$C(n+1,d,j)=\sum_{k=0}^{n+1}(-1)^k\binom{n+1}{k}
\binom{j-1-k(d-1)}{n}.$}$$
Here $\binom{p}{q}=0$ if $q<0$ or $p-q<0$.

On the other hand, we have by a well-known result of Griffiths [12]
$$C(n+1,d,di)=\dim \Gr_F^{n-i}H^{n-1}_{\prim}(Z,\C)\q\h{for}
\,\,\, i\in [1,n].$$
Using these, we can calculate the Hodge numbers of the limit mixed
Hodge structure in Corollary~5.

\ms
(iii) With the notation of Theorem~5 we have the weight spectral
sequence for the vanishing cycles
$$E_1^{-k,j+k}=H^j(X_0,\Gr^W_{k+n-1}\ph_f(\Q_X[n]))\Rightarrow
H^j(X_0,\ph_f(\Q_X[n])),$$
degenerating at $E_2$, and $E_1^{-i,j+i}$ is given by the
direct sum of
$$\bigoplus_{l\ge\max(-i,1)}\Bigl(\bigoplus_{|I|=i+2l+1}
H^{j+n-|I|}(Y_I)(-i-l)\Bigr)$$ and
$$\bigoplus_{l\ge\max(-i,0)}\Bigl(\bigoplus_{|I|=i+2l+2}
H^{j+n-|I|-1}(Y'_I)(-i-l-1)\Bigr).$$

\ms
(iv) Let $i:X\into Y$ be a closed immersion of smooth varieties.
Assume $\codim_YX=1$.
Set $n=\dim Y$.
Since $i^*\Q_Y=\Q_X$ and $i^!\Q_Y=\Q_X(-1)[-2]$, the adjunction
for $i$ implies
$$\Ext^1(\Q_Y[n],\Q_X[n-1])=\Ext^1(\Q_X[n-1],\Q_Y(1)[n])=
\Hom(\Q_X,\Q_X).$$

\ms
(v) For a quasi-semistable degeneration, it would be possible to
get a semistable model by repeating the blow-ups along the center
$\{t=y_k=0\}$ in the notation of (3.1.1).
Indeed, the strict transform of $\{y_1\cdots y_k=y_nt\}$ becomes
$\{y_1\cdots y_{k-1}=y_nt\}$ or $\{y_1\cdots y_k=y_n\}$
by replacing $t$ with $ty_k$ or $y_k$ with $ty_k$, see [13] for
the case $r=2$.
If $r=2$ we can prove Theorem~(3.2) using this blow-up,
see [20].

\section{Application}

\nin
In this section we prove Theorem~6 which improves [15], Cor.~6.4.

\ms\nin
{\bf 4.1.~Proof of Theorem~6.}
By assumption $L=L_1\otimes L_2$ with $L_1,L_2$ very ample.
Moreover $L_2$ is a tensor product of two very ample line bundles
if $Y=\P^n$ with $n$ even.
Let $Y_i\in|L_i|\,(i=0,1,2)$ be general smooth members where
$L_0:=L$.
We may assume that their union is a divisor with normal crossings.
Let $g_i\in\Gamma(X,L_i)$ defining $Y_i$ for $i=0,1,2$.
Applying the construction of Theorem~5, we get a smoothing
$f:X\to\C$ of $Y_1\cup Y_2$ such that
$$X_t:=f^{-1}(t)=\{g_1g_2=tg_0\}\subset X.$$
In particular, $X_0=Y_1\cup Y_2$ and $X_1=Y_0$.
The singular locus $\S$ of $X$ is contained in $X_0=\{t=0\}$ and
coincides with the intersection of $Y_0,Y_1,Y_2$, see (3.1.2).
Set $Z=Y_1\cap Y_2$.
Note that $\dim Y_i=n-1$, $\dim Z=n-2$ and $\dim\S=n-3$.
By Theorem~5 we have the weight spectral sequence
$$E_1^{-k,j+k}=H^j(X_0,\Gr^W_{k+n-1}\psi_{f,1}\IC_X\Q)\Rightarrow
H^{j+n-1}(X_{\infty},\Q),$$ degenerating at $E_2$, and
its $E_1$-terms are given as follows:
$$E_1^{-k,j+k}=\begin{cases}H^{j+n-2}(Z)(-1)& \text{if
$\,k=1$,}\\ \mopl_{i=1,2\,}H^{j+n-1}(Y_i)\oplus H^{j+n-3}(\S)(-1)
&\text{if $\,k=0$,}\\ H^{j+n-2}(Z)&\text{if $\,k=-1$,}\\ \,0&
\text{otherwise.}\end{cases}$$
(It is also possible to get this spectral sequence by the same
argument as in [20], 4.4 using the blow-up along
$Y_0\cap Y_2$ studied in [13].)
The differential $d_1^{p,q}$ is expressed by the restriction and
Gysin morphisms associated to the inclusions $\S\to Z$ and
$Z\to Y_i$ up to constant multiples,
see see Remark~(3.5)(iv) below. Set
$$\aligned\a_r^{j,k}&=\dim E_r^{-k,j+k},\\
\g^j(V)&=\dim H^{j+\dim V}(V)\,\,\,\h{for}\,\,V=Y,Y_i,Z,\S.
\endaligned$$
There are symmetries
$$\a_1^{j,k}=\a_1^{-j,k}=\a_1^{j,-k}=\a_1^{-j,-k},\q
\g^{-j}(V)=\g^j(V).$$
We have to show
$$\g^0(Y_0)=\msum_{|k|\le 1}\,\a_2^{0,k}>\g^{-1}(Y),
\leqno(4.1.1)$$
where the first equality follows from the $E_2$-degeneration of
the above spectral sequence.
We will show (4.1.1) by induction on $n$.

The weak Lefschetz theorem implies
$$\aligned\g^0(\S)\ge\g^{-1}(Z)=\g^{-2}(Y_i)=\g^{-3}(Y),\\
\g^0(Z)\ge\g^{-1}(Y_i)=\g^{-2}(Y),\\
\g^0(Y_i)\ge\g^{-1}(Y).\endaligned\leqno(4.1.2)$$
From the above spectral sequence we get then
$$\aligned&\a_2^{0,0}-\g^{-1}(Y)
\ge\a_1^{0,0}-\a_1^{-1,1}-\a_1^{1,-1}-\g^{-1}(Y)\\
&=\g^0(Y_1)-\g^{-1}(Y)+\g^0(Y_2)-\g^{-3}(Y)+\g^0(\S)-\g^{-3}(Y),
\endaligned\leqno(4.1.3)$$
using $\a_1^{-1,1}=\a_1^{1,-1}=\g^{-1}(Z)=\g^{-3}(Y)$.
We have moreover
$$\a_2^{0,-1}\ge\g^0(Z)-\g^{-2}(Y)=\g^0(Z)-\g^{-1}(Y_1).
\leqno(4.1.4)$$
Indeed, by the weak Lefschetz theorem, the image of $H^{n-2}(Y_i)
\,(i=1,2)$ in $H^{n-2}(Z)$ under the differential $d_1$ of the
spectral sequence is contained in the image of the canonical
injection $H^{n-2}(Y)\into H^{n-2}(Z)$ by the restriction
morphism.
Moreover, the image of $H^{n-4}(\S)(-1)$ in $H^{n-2}(Z)$ is also
contained in the image of the above canonical injection since it
is contained in the image of the action of $c_1(L_0)$ using the
bijectivity of the restriction morphism $H^{n-4}(Z)\to H^{n-4}(\S)$
in the case $n\ge 4$.
(Note that the action of $c_1(L_0)$ is compatible with the
restriction morphism.)
So (4.1.4) follows.

Thus, in order to prove (4.1.1) using (4.1.2--4), it is enough
to show either
$$\g^0(\S)>\g^{-1}(Z)=\g^{-3}(Y)\q\h{or}\q
\g^0(Z)-\g^{-1}(Y_1)>0.\leqno(4.1.5)$$
Using the first condition of (4.1.5), the assertion is reduced
to the case where $n$, $Y$ and $Y_0$ are respectively replaced by
$n-2$, $Z$ and $\S$.
So we get the assertion in the $n$ odd case by induction on $n$.
Indeed, in case $n=3$ we have $\#|\S|\ge 2$ since it coincides
with the intersection number $Y_0\cdot Y_1\cdot Y_2$ where $Y_0$
can be replaced by $Y_1+Y_2$ as algebraic cycles and $Y_1,Y_2$
are very ample by assumption.

We may now assume $n$ even.
In the case where $Y=\P^n$ with $n$ even so that $L$ is a tensor
product of three very ample line bundles, the assertion then
follows from the $n$ odd case using the last condition of (4.1.5)
and replacing $Y$ with $Y_1$.
In the other case, the assertion is reduced finally to the case
$n=2$ by induction on $n$ using the first condition of (4.1.5)
repeatedly where $n$, $Y$ and $Y_0$ are respectively replaced
by $n-2$, $Z$ and $\S$ in each inductive step.
We have then $\S=\emptyset$, $\dim Z=0$ and $Y_i$ is a connected
curve where the image of $H^0(Y_i)\to H^0(Z)$ is the diagonal
for $i=1,2$ and its cokernel is nonzero if $\#|Z|>1$.
The last condition is satisfied after replacing $Y$ with $Z$
repeatedly if the intersection number $Y_1^{n/2}\cdot Y_2^{n/2}$
is bigger than 1 for the original $Y_1,Y_2\subset Y$.
So the assertion is reduced to the following lemma (which would be
known to specialists at least in the case $D_i=D_j$ for any $i,j$).

\ms\nin
{\bf Lemma~4.2.} {\it Let $Y$ be a smooth complex projective
variety of dimension $n\ge 2$, and $D_i$ be very ample divisors
for $i=1,\dots,n$.
Assume the intersection number $D_1\cdots D_n$ is $1$.
Then $Y=\P^n$ and the $D_i$ are hyperplanes.}

\ms\nin
{\it Proof.}
Note first that the assertion is easy if $D_i\in|L|$ for some very
ample line bundle $L$ independent of $i$.
Indeed, take a linear subsystem generated by general $n+1$
hyperplane sections and defining a morphism $Y\to\P^n$.
Then the assumption implies that its fiber is one point over any
$s\in \P^n$ (taking $n$ hyperplanes in $\P^n$ whose intersection
is $s$).
So the assertion follows.

We prove the general case by induction on $n\ge 1$.
If $n=1$, the intersection number is interpreted as the degree
of a zero-cycle.
Then the assertion is well known.
Assume $n\ge 2$.
We first show that the $D_i$ are isomorphic to $\P^{n-1}$.
For this it is enough to show that the $D_i$ are smooth using the
inductive hypothesis.
But it is well known (and is easy to show) that the intersection
number cannot be 1 if some $D_i$ has a singular point.
Note that $D_i$ may be replaced by any member of the linear
system $|\OO_Y(D_i)|$ and the $D_i$ are very ample.

We now consider a Lefschetz pencil $f:X\to \P^1$ associated to the
very ample line bundle $\OO_Y(D_n)$.
Its fibers are {\it all} isomorphic to $\P^{n-1}$ by the
above argument.
So we get
$$\h{$h^{2i}(X)=2$ if $i=1,\dots,n-1$, and $h^j(X)=0$ if $j$ is
odd}.$$
Since the morphism $\rho:X\to Y$ is the blow-up along a smooth
center $C$ which has codimension 2, we have
$$h^{2i}(X)=h^{2i}(Y)+h^{2i-2}(C)\q\h{for any}\,\,i,$$
and $h^{2i}(Y)$ and $h^{2i-2}(C)$ are nonzero for $i=1,\dots,n-1$
since $Y,C$ are projective. We get thus
$$\h{$h^{2i}(Y)=1$ if $i=1,\dots,n-1$, and $h^j(Y)=0$ if $j$ is
odd}.$$
Moreover, $H^2(Y,\Z)$ is torsion-free since so is $H^2(X,\Z)$
and there is a long exact sequence
$$0\to H^2(Y,\Z)\to H^2(X,\Z)\to H^0(C,\Z)(-1)\to.$$
(The last sequence is induced by the truncation $\tau$ on the
direct image $\R\rho_*\Z_X$ since $R^i\rho_*\Z_X$ is $\Z_Y$ if
$i=0$, $\Z_C(-1)$ if $i=2$, and 0 otherwise.)
So we get ${\rm Pic}(Y)=\Z[D]$ where $D$ is an ample divisor,
and $[D_i]=m_i[D]$ for positive integers $m_i$.
Then the assumption implies $m_i=1$ for any $i$,
and the assertion is reduced to the first case.

\ms\nin
{\bf Remarks~4.3.}
(i) The equivalence between the non-vanishing and the
non-surjectivity in Theorem~6 follows from the Picard-Lefschetz
formula [17] and the global invariant cycle theorem [5],
see [15] (and [16] for a complex analytic argument).

\ms
(ii) If $Y=\P^n$ with $n$ even and $L=\OO_{\P^n}(2)$, then
$H^{n-1}(Y_s)=0$ for a smooth member $Y_s$ of $|L|$ and we have
the surjectivity of $i_s^*$.
So Theorem~6 is optimal (i.e.\ we cannot take $k=2$) in the case
$Y$ is projective space of even dimension.
Among the other cases where Theorem~6 holds with $k=2$, there
are some cases where we cannot take $k=1$.
For example, if $Y$ is a projective hypersurface of degree
$d\le 2$ and $L=\OO_Y(1)$, then we have the surjectivity of $i_s^*$
except for the case where $\dim Y$ is odd and $d=2$.
Recently N. Fakhruddin informed us that there are less trivial
examples.
For instance, if $Y$ is a ruled surface over $\P^1$, then we have
a very ample divisor which is a sum of a section $C$ and the
pull-back of a very ample divisor on $\P^1$.
In this case any member of the linear system $|L|$ is a section
(since its intersection number with any fiber is 1) and
we have the surjectivity of $i_s^*$ for $k=1$.

\ms
(iii) The non-surjectivity of the restriction morphism $i_s^*$
in Theorem~6 is equivalent to the condition that
$H^n(Y\setminus Y_s)$ is not pure of weight $n$ by Poincar\'e
duality and the weight spectral sequence [5].
This implies examples of smooth affine varieties $Y\setminus Y_s$
whose cohomology is not pure and which are formal in the sense of
[8], see [19], Cor.~10.5, a).

\ms
(iv) Let $\pi:\Y\to S:=(\P^N)^{\vee}$ denote the universal family of
hyperplane sections as in a remark after Theorem~6.
By [3] we have the decomposition
$$\R\pi_*(\Q_{\Y}[\dim\Y])\cong\mopl_{i,j}E_j^i[-i],$$
where $E_j^i$ denotes the direct sum of the direct factors of
${}^p\!R^i\pi_*(\Q_{\Y}[\dim\Y])$ having strict support of
codimension $j$.
In [4] it is proved that
$$E_j^i=0\q\h{for}\,\,\,ij\ne 0.$$
Moreover, N.~Fakhruddin proved in Appendix of loc.~cit.
$$E_j^0=0\,\,\,\h{for $j\ge 1$ if $L=L'{}^{\otimes d}$ with
$L'$ ample and $d\gg 1$.}$$
Here it is sufficient to assume $d\ge 2n-1$ if $L'$ is very ample.
Indeed, let $\I_Z\subset\OO_Y$ be the reduced ideal sheaf of
$Z=\{z_1,\dots,z_m\}\subset Y$.
By the same argument as in loc.~cit., the assertion can be
reduced to the surjection
$$H^0(Y,L'^{\otimes d})\onto H^0(Y,L'^{\otimes d}\otimes_{\OO_Y}
(\OO_Y/\I_Z^2))\q\h{for}\,\,\,d\ge 2m-1.$$
For the last surjection, take $H_{i,j}\in|L|$ for $i\in[1,m]$,
$j\in[0,n]$ such that $z_i\notin H_{i,0}$ and the $H_{i,j}$ for
$j\in[1,n]$ contain $z_i$ and form a divisor with normal crossings
at $z_i$.
Take further $H'_i\in|L^{m-1}|$ such that $z_i\notin H'_i$ and
$z_j\in H'_i$ for $j\ne i$.
Then the surjection follows by using $H_{i,j}+2H'_i\in|L^{2m-1}|$
for $i\in[1,m]$, $j\in[0,n]$.
Note that Theorem~6 implies (assuming $Y\ne\P^n$)
$$E_1^0=0\,\,\,\h{if $L=L'^{\otimes d}$ with $L'$ very ample
and $d\ge 2$.}$$

\ms
(v) Let $Y_s$ be a smooth fiber of $\pi$ in Remark~(iv).
Set $n=\dim Y$.
Using the Picard-Lefschetz formula, we have the orthogonal
decomposition
$$H^{n-1}(Y_s)=H^{n-1}(Y_s)^{\inv}\oplus H^{n-1}(Y_s)^{\van},$$
where $H^{n-1}(Y_s)^{\inv}=\Im(H^{n-1}(Y)\to H^{n-1}(Y_s))$ and
$H^{n-1}(Y_s)^{\van}$ is generated by the vanishing cycles
(this is closely related to the hard Lefschetz theorem,
see [7], [16].)

We have
$$H^{n-1}(Y_s)^{\van}\subset H^{n-1}_{\prim}(Y_s),$$
where $H^{n-1}_{\prim}(Y_s)$ denotes the primitive part.
Indeed, the non-primitive part is contained in the invariant part
and the assertion follows by taking their orthogonal complements.

Let $f:X\to\P^1$ be a Lefschetz pencil where $\rho:X\to Y$ is the
blow-up along a smooth center $Z$ which is the intersection
of two general hyperplane sections.
We have the perverse Leray spectral sequence
$$E_2^{i.j}=H^i(\P^1,{}^p\!R^jf_*(\Q_X[n]))\Rightarrow
H^{i+j+n}(X),$$
degenerating at $E_2$ by the decomposition theorem [3].
Note that
$$\aligned E_2^{-1,j+1}&=\Im(\i^*_s:H^{j+n}(X)\to H^{j+n}(Y_s)),
\\ E_2^{1,j-1}&=\Im((\i_s)_*:H^{j+n-2}(Y_s)(-1)\to H^{j+n}(X)),
\endaligned$$ where $\i_s:Y_s\into X$ is the inclusion.
Moreover we have the decomposition
$$H^{j+n}(X)=H^{j+n}(Y)\oplus H^{j+n-2}(Z)(-1),$$
so that $\i^*_s,(\i_s)_*$ are expressed by using $i^*_s,(i_s)_*,
i'^*_s,(i'_s)_*$ where $i_s:Y_s\into Y$, $i'_s:Z\into Y_s$
are the inclusions, see [15].
Using the weak Lefschetz theorem, we get then a canonical
isomorphism as Hodge structures
$$E_2^{0,0}=H^n_{\prim}(Y)\oplus H^{n-2}(Z)^{\van}(-1).$$

Let $D\subset S$ denote the discriminant of $\pi$ in Remark~(iv).
In the case $\codim_SD>1$, $^p\!R^0f_*(\Q_X[n])$ is constant so that
$E_2^{0,0}=0$, and we get by using the above formula
$$H^n_{\prim}(Y)=H^{n-1}(Y_s)^{\van}=H^{n-2}(Z)^{\van}=0.$$
Note that its converse is also true.
Indeed, the equivalence can be shown by using the following
(see [16], 3.5.3)
$$\chi(Y)-2\chi(Y_s)+\chi(Z)=(-1)^nr,$$
where $r$ is the number of the critical values of the Lefschetz
pencil, i.e. the intersection number of $D$ with a general line.
The weak Lefschetz theorem implies that the left-hand side
coincides up to a sign with a sum of three non-negative numbers
$$\bigl(b_n(Y)-b_{n-2}(Y)\bigr)+ 2\bigl(b_{n-1}(Y_s)-b_{n-1}(Y)
\bigr)+\bigl(b_{n-2}(Z)-b_{n-2}(Y_s)\bigr).$$

\end{document}